\documentclass[12pt]{article}
\usepackage{mathrsfs}
\usepackage{pifont}

\usepackage{amsfonts}
\usepackage{latexsym}
\usepackage{amsmath}
\usepackage{amssymb}
\usepackage{color}

\def\a{\alpha}
\def\ep{\varepsilon}
\def\p{\varphi}
\def\d{\delta}

\def\s{\sigma}

\def\g{\gamma}
\def\l{\lambda}

\title{On the comparison theorem for multidimensional SDEs with jumps}

\date{June 7, 2010}

\author{ Xuehong Zhu}

\begin{document}

\maketitle

\begin{abstract}
In this note, we give a necessary and sufficient condition under which the comparison theorem holds for multidimensional stochastic
differential equations (SDEs) with jumps and for matrix-valued SDEs with jumps.\\
\par  $\textit{Keywords: Vibility property; Viscosity solution; Comparison theorem.}$
\end{abstract}



\section{Introduction}\label{sec:intro}
\qquad The comparison theorem for real-valued SDEs turns out to be
one of the classic results of this theory. We can refer the reader
to [1], [2], [6], [7] and so on. It allows to compare the solutions
of two real-valued SDEs whenever we can compare the drift and the
diffusion coefficients. Thus they are all sufficient conditions.

Until in [4], Peng and Zhu originally studied comparison theorem of
1-dimensional SDEs with jumps through the viability theory (see Peng
and Zhu [5]) and got the necessary and sufficient condition. In the
manuscirpt of [3], Hu and Peng studied the multidimensional
situation without jumps applying the viability theory.

The objective of this paper is to give a necessary and sufficient
condition under which the comparison theorem holds for
multidimensional SDEs with jumps. For this, we still apply the
stochastic viability property (SVP) for SDEs with jumps studied in
[5] and combine the technique used in [3].

The paper is organized as follows: in the next section, we recall
the viability criteria for SDEs with jumps; in section 3, we study
the comparison theorem for multidimensional SDEs with jumps and for
matrix-valued SDEs with jumps.

\section{A characterization for SDEs with jumps under state constraint}\label{sec:intro}
\qquad Let $(\Omega,{\cal{F}},P,({\cal{F}}_{t})_{t\geq 0})$ be a
complete stochastic basis such that $\mathcal{F}_{0}$ contains all
$P$-null elements of ${\cal{F}}$, and
$\mathcal{F}_{t^{+}}:=\cap_{\ep>0}\mathcal{F}_{t+\ep}=\mathcal{F}_{t},t\geq
0$, and $\mathcal{F}=\mathcal{F}_{T}$, and suppose that the
filtration is generated by the following
two mutually independent processes:\\
(i) a $d$-dimensional standard Brownian motion $(W_{t})_{0\leq
t\leq T}$, and\\
(ii) a stationary Poisson random measure $N$ on $(0,T]\times E$,
where $E\subset R^{l}\setminus\{0\}$, $E$ is equipped with its Borel
field ${\mathscr{B}}_{E}$, with compensator $\hat{N}(dtde)=dtn(de)$,
such that $n(E)<\infty$, and $\{\tilde{N}((0,t]\times
A)=(N-\hat{N})((0,t]\times A)\}_{0<t\leq T}$ is an
$\mathcal{F}_{t}$-martingale, for each $A\in \mathscr{B}_{E}$.

By $T>0$ we denote the finite real time horizon.

We consider a jump diffusion process as follows:
$$
X^{t,x}_{s}=x+\int^{s}_{t}b(r,X^{t,x}_{r})dr+\int^{s}_{t}\sigma(r,X^{t,x}_{r})dW_{r}\\
+\int^{s}_{t}\int_{E}\gamma(r,X^{t,x}_{r-},e)\tilde{N}(drde),s\in
[t,T], \eqno(2.1)
$$
where
$$
\begin{array}{l}
b:[0,\infty)\times R^{m} \rightarrow R^{m},
\gamma:[0,\infty)\times R^{m} \times R^{l} \rightarrow R^{m},\\
\sigma=\{\sigma_{\alpha}^{i}\}:[0,\infty)\times R^{m} \rightarrow
R^{m\times d},i=1,2,...,m,\a=1,2,...,d.
\end{array}
$$

{\bf Definition 2.1.} The SDE (2.1) enjoys the stochastic viability
property (SVP in short) in a given closed set  $K\subset R^{m}$ if
and only if: for any fixed time interval $[0,T]$, for each $(t,x)\in
[0,T]\times K$, there exists a probability space
 $(\Omega,{\cal{F}},P)$, a $d-$dimensional Brownian motion $W$, a
 stationary Poisson process $N$,  such that
 $$
 X^{t,x}_{s}\in K ,\mbox{ \ } \forall \mbox{ \ }s  \in [t,T] \mbox{ \
 } \mbox{ \ }P-a.s..
 $$

We assume that, there exists a sufficiently large constant $\mu>0$
and a function $\rho: R^{l}\rightarrow R_{+}$ with
$$
\int_{E}\rho^{2}(e)n(de)<\infty,
$$
such that

(A1)$b,\sigma,\gamma\mbox{ \ }\mbox{are continuous in}\mbox{ \ }
(t,x),$

(A2) for all $x,x'\in R^{m}$, $t\in [0,+\infty)$
$$
\begin{array}{ll}
| b(t,x)-b(t,x')|+|\s(t,x)-\s(t,x')|\leq \mu|x-x'|,\\
|b(t,x)|+|\s(t,x)|\leq \mu(1+|x|),\\
|\g(t,x,e)-\g(t,x',e)|\leq \rho(e)|x-x'|,\forall e\in E,\\
|\g(t,x,e)|\leq \rho(e)(1+|x|),\forall e\in E.
\end{array}
$$
Here $\langle\cdot\rangle$ and $|\cdot|$ denote, respectively, the
Euclidian scalar product and norm. Obviously under the above
assumptions there exists a unique strong solution to SDE (2.1). $C$
is a constant such that
$$
C\geq 1+2\mu+\mu^{2}+\int_{E}\rho^{2}(e)n(de).
$$

We denote by $C_{2}([0,T]\times R^{m})$ (resp,.
$C^{1,2}_{2}([0,T]\times R^{m})$) the set of all functions in
$C([0,T]\times R^{m})$ (resp., $C^{1,2}([0,T]\times R^{m})$) with
quadratic growth in $x$. In fact, the SVP in $K$ is related to the
following PDE:
$$
 \left\{
\begin{array}{l}
\mathscr{L}u(t,x)+\mathscr{B}u(t,x)-Cu(t,x)+d^{2}_{K}(x)=0,\mbox{ \
}(t,x)\in (0,T)\times R^{m},\\
u(T,x)=d^{2}_{K}(x),
\end{array}
\right. \eqno{(2.2)}
$$
where we denote, for $\p\in C_{2}^{1,2}([0,T]\times R^{m})$ ,
$$
\mathscr{L}\p(t,x):=\frac{\partial \p(t,x)}{\partial t}+\langle
D\p(t,x),b(t,x)\rangle +\frac{1}{2}tr[D^{2}\p(t,x)\s\s^{T}(t,x)],
$$
$$
\mathscr{B}\p(t,x):=\int_{E}[\p(t,x+\g(t,x,e))-\p(t,x)-\langle
D\p(t,x),\g(t,x,e)\rangle]n(de).
$$

{\bf Definition 2.2.} We say a function $u\in C_{2}([0,T]\times
R^{m})$ is a viscosity supersolution (resp., subsolution) of (2.2)
if, $u(T,x)\geq d^{2}_{K}(x)$ (resp., $u(T,x)\leq d^{2}_{K}(x)$) and
for any $\p\in C_{2}^{1,2}([0,T]\times R^{m})$ and any point
$(t,x)\in [0,T]\times R^{m}$ at which $u-\p$ attains its minimum
(resp., maximum),
$$
\mathscr{L}\p(t,x)+\mathscr{B}\p(t,x)-C\p(t,x)+d^{2}_{K}(x)\leq
0,\mbox{ \ (resp.,  }\geq 0).
$$
u is called a viscosity solution if it is both viscosity
supersolution and subsolution.

Now let us recall the characterization of SVP of SDE (2.1) in $K$
(see [5]):

{\bf Lemma 2.3.}We assume (A1) and (A2). Then the following claims
are
equivalent:\\
{\rm(i)}SDE (2.1) enjoys the SVP in $K$; \\
{\rm(ii)}$d^{2}_{K}(\cdot)$ is a viscosity supersolution of PDE
(2.2).

\section{Comparison theorem for SDEs with jumps}\label{sec:intro}
\subsection{multidimensional SDEs}\label{sec:intro}
\qquad Let $\mathcal{S}^{2}_{[0,T]}$ denote the set of
$\mathcal{F}_{t}$-adapted c$\grave{a}$dl$\grave{a}$g $m$-dimensional
processes\\ $\{X_{t},0\leq t\leq T\}$ which are such that
$$
\|X\|_{\mathcal{S}^{2}_{[0,T]}}:=(E[\sup_{0\leq t\leq
T}|X_{t}|^{2}])^{\frac{1}{2}}<\infty.
$$

Consider the following two SDEs: $i=1,2,$
$$
X^{i}_{s}=x^{i}+\int_{t}^{s}b^{i}(r,X^{i}_{r})ds+\int_{t}^{s}\s^{i}(r,X^{i}_{r})dW_{r}+\int_{t}^{s}\int_{E}\g^{i}(r,X^{i}_{r-},e)\tilde{N}(drde),
\eqno{(3.1)}
$$
where $(b^{i},\s^{i},\g^{i}),i=1,2$ satisfy (A1) and (A2), and
$x^{1},x^{2}\in R^{m}$. The objective of this section is to study
when the comparison theorem holds for two SDEs with jumps of type
(3.1). We will find that the comparison theorem can be transformed
to a viability problem in $R^{m}_{+}\times R^{m}$ of
$(X^{1}-X^{2},X^{2})$.

{\bf Theorem 3.1.} Suppose that $(b^{i},\s^{i},\g^{i}),i=1,2$
satisfy (A1) and (A2). Then the following are equivalent:

{\rm(i)} For any $t\in [0,T]$, $x^{1},x^{2}\in R^{m}$ such that
$x^{1}\geq x^{2}$, the unique adapted solutions $X^{1}$ and $X^{2}$
in $\mathcal{S}^{2}_{[t,T]}$ to the SDE (3.1) over time interval
$[t,T]$ satisfy:
$$
X^{1}_{s}\geq X^{2}_{s},s\in[t,T],P-a.s.;
$$

{\rm(ii)} $\s^{1}\equiv \s^{2}$, and for any $t\in
[0,T],k=1,2,...,m,$
$$
\left\{\begin{array}{ll}
{\rm(a)}& \s^{1}_{k} \mbox{ \ depends only on \ }x_{k},\\
{\rm(b)}&\forall x'\in R^{m}, \forall x\geq
0,x_{k}+\g^{1}_{k}(t,x+x',e)-\g^{2}_{k}(t,x',e)\geq 0,n(de)-a.s.,\\
{\rm(c)}& \mbox{for all \ }x',\d^{k}x\in R^{m}, \mbox{ \ such that
\ }\d^{k}x\geq0,(\d^{k}x)_{k}=0,\\
& b^{1}_{k}(t,\d^{k}x+x')-\displaystyle\int_{E}\g^{1}_{k}
(t,\d^{k}x+x',e)n(de)\geq
b^{2}_{k}(t,x')-\displaystyle\int_{E}\g^{2}_{k} (t,x',e)n(de).
\end{array}\right.
$$

{\bf Proof: }Set $\bar{X}_{s}=(X^{1}_{s}-X^{2}_{s},X^{2}_{s})$, then
(i) is equivalent to the following:

For any $t\in [0,T]$, $\forall \bar{x}=(x^{1},x^{2})$ such that
$x^{1}\geq 0$, the unique solution $\bar{X}$ to the following SDE
over time interval $[t,T]$:
$$
\bar{X}_{s}=\bar{x}+\int_{t}^{s}\bar{b}(r,\bar{X}_{r})ds+\int_{t}^{s}\bar{\s}(r,\bar{X}_{r})dW_{r}+\int_{t}^{s}\int_{E}\bar{\g}(r,\bar{X}_{r-},e)\tilde{N}(drde),
\eqno{(3.2)}
$$
satisfies $\bar{X}^{1}_{s}\geq 0,s\in [t,T],P-a.s.$, where for
$\bar{x}=(\bar{x}^{1},\bar{x}^{2})$,
$$
\begin{array}{l}
\bar{b}(s,\bar{x})=(b^{1}(s,\bar{x}^{1}+\bar{x}^{2})-b^{2}(s,\bar{x}^{2}),b^{2}(s,\bar{x}^{2})),\\
\bar{\s}(s,\bar{x})=(\s^{1}(s,\bar{x}^{1}+\bar{x}^{2})-\s^{2}(s,\bar{x}^{2}),\s^{2}(s,\bar{x}^{2})),\\
\bar{\g}(s,\bar{x},e)=(\g^{1}(s,\bar{x}^{1}+\bar{x}^{2},e)-\g^{2}(s,\bar{x}^{2},e),\g^{2}(s,\bar{x}^{2},e)).
\end{array}
$$

So we can apply Lemma 2.3 to SDE (3.2) and the convex closed set
$K:=R^{m}_{+}\times R^{m}$, i.e., (i) is equivalent to that
$d^{2}_{K}(\cdot)$ is a viscosity supersolution of PDE (2.2).

we can see that $\forall x=(x^1,x^2)\in R^{2m}$,
$$
\Pi_{K}(x)=\left(\begin{array}{c}
(x^{1})^{+}\\\\
x^{2}
\end{array}\right),
x-\Pi_{K}(x)=\left(\begin{array}{c}
-(x^{1})^{-}\\\\
0
\end{array}\right).
$$
So
$$
d^{2}_{K}(x)=|(x^{1})^{-}|^{2}=\sum\limits_{k=1}^{m}I_{\{x^{1}_{k}<0\}}|x^{1}_{k}|^{2}.
$$
And
$$
(D^{2}d^{2}_{K})(x)\left\{\begin{array}{ll} =0_{2m\times 2m},&
\mbox{when \ }x\in K^\circ,\\
\mbox{does not exist},& \mbox{when \ }x\in \partial K,\\
=(a_{ij})_{2m\times 2m},& \mbox{when \ }x\in R^{2m}\backslash K,
\end{array}
\right.
$$
where
$$
a_{ij}=0,\mbox{ \ when \ }i\neq j,\mbox{ \ } a_{ii}=\left\{
\begin{array}{ll}
0,& m<i\leq 2m,\\
0,& 1\leq i\leq m,x^{1}_{i}\geq 0,\\
2,& 1\leq i\leq m,x^{1}_{i}< 0.
\end{array}
\right.
$$

From the above analysis and the Lipschitz condition of $b^{1}$
w.r.t. $x$, we can easily check that: $d^{2}_{K}(\cdot)$ is a
viscosity supersolution of PDE (2.2) if and only if, \\\\
$\mbox{(ii)}'$ $\forall t\in [0,T], \forall (x,x')\in R^{m}\times
R^{m}$,
$$
\begin{array}{ll}
&-2\langle
x^{-},b^{1}(t,x^++x')-b^{2}(t,x')\rangle+\sum\limits_{k=1}^{m}I_{\{x_{k}<0\}}|\s^{1}_{k}(t,x+x')-\s^{2}_{k}(t,x')|^{2}\\
&+\sum\limits_{k=1}^{m}I_{\{x_{k}<0\}}\displaystyle\int_{E}[|(x_{k}+\g_{k}^{1}(t,x+x',e)-\g_{k}^{2}(t,x',e))^{-}|^{2}\\
&\mbox{ \ \ \ \ \  \ \ \ \ \ \ \ \ \ \ \ \ \ \ \ }-|x_{k}|^{2}-2x_{k}(\g_{k}^{1}(t,x+x',e)-\g_{k}^{2}(t,x',e))]n(de)\\
&+\sum\limits_{k=1}^{m}I_{\{x_{k}\geq
0\}}\displaystyle\int_{E}|(x_{k}+\g_{k}^{1}(t,x+x',e)-\g_{k}^{2}(t,x',e))^{-}|^{2}n(de)\\
\leq & C^*|x^{-}|^{2},
\end{array}
$$
where $C^*\geq 4\mu+\mu^{2}+\int_{E}\rho^{2}(e)n(de)$ is a constant
which does not depend on $t,x,x'$. Then the left thing we need to do
is to prove: $\mbox{(ii)}\Leftrightarrow \mbox{(ii)}'$.

$\mbox{(ii)}'\Rightarrow\mbox{(ii)}$: If we pick $x\geq 0$, we can
immediately get (b) in (ii) from (ii)$'$.

Pick $x<0$, by (ii)$'$ we have
$$
\sum\limits_{k=1}^{m}2x_{k}
[b^{1}_{k}(t,x')-b^{2}_{k}(t,x')]+|\s^{1}(t,x+x')-\s^{2}(t,x')|^{2}\leq
C^{*}|x|^{2}.
$$
Let $x$ tend to $0_{-}$, we get that $\s^{1}\equiv \s^{2}$.

And $\forall \d^{k}x\in R^{m}$, such that,
$\d^{k}x\geq0,(\d^{k}x)_{k}=0$. Pick $x=\d^{k}x-\ep e_{k},\ep>0$.
From (ii)$'$ we get:
$$
-2\ep
[b^{1}_{k}(t,\d^{k}x+x')-b^{2}_{k}(t,x')]+|\s^{1}_{k}(t,x+x')-\s^{1}_{k}(t,x')|^{2}\leq
C^{*}\ep^{2}.
$$
Let $\ep$ tend to 0, we have
$$
\s^{1}_{k}(t,\d^{k}x+x')=\s^{1}_{k}(t,x').
$$
We deduce quickly that $\s^{1}_{k}$ depends only on $x_{k}$.

With $x=\d^{k}x-\ep e_{k},\ep>0$ again, from (ii)$'$ we can also get
$$
-2\ep
[b^{1}_{k}(t,\d^{k}x+x')-b^{2}_{k}(t,x')]+\int_{E}[-\ep^{2}+2\ep(\g^{1}_{k}(t,x+x',e)-\g^{2}_{k}(t,x',e))]n(de)\leq
C^{*}\ep^{2}.
$$
Dividing by $-2\ep$ and letting $\ep$ tend to 0, we have
$$
b^{1}_{k}(t,\d^{k}x+x')-\int_{E}\g^{1}_{k} (t,\d^{k}x+x',e)n(de)\geq
b^{2}_{k}(t,x')-\int_{E}\g^{2}_{k} (t,x',e)n(de).
$$

(ii)$\Rightarrow$(ii)$'$. For $x\geq 0$, from (b) in (ii), we know
that (ii)$'$ holds true. If there exist some $1\leq k\leq n$, such
that $x_{k}<0$, then from (ii) we have
$$
\begin{array}{ll}
&-2\langle
x^{-},b^{1}(t,x^{+}+x')-b^{2}(t,x')\rangle+\sum\limits_{k=1}^{m}I_{\{x_{k}<0\}}|\s^{1}_{k}(t,x+x')-\s^{1}_{k}(t,x')|^{2}\\
&+\sum\limits_{k=1}^{m}I_{\{x_{k}<0\}}\displaystyle\int_{E}[|(x_{k}+\g^{1}_{k}(t,x+x',e)-\g^{2}_{k}(t,x',e))^{-}|^{2}\\
&\mbox{ \ \ \ \ \ \ \ \ \ \ \ \ \  \ \ \ \ \ \ \ }-|x_{k}|^{2}-2x_{k}(\g^{1}_{k}(t,x+x',e)-\g^{2}_{k}(t,x',e))]n(de)\\
&+\sum\limits_{k=1}^{m}I_{\{x_{k}\geq
0\}}\displaystyle\int_{E}|(x_{k}+\g^{1}_{k}(t,x+x',e)-\g^{2}_{k}(t,x',e))^{-}|^{2}n(de)\\
\leq &
\sum\limits_{k=1}^{m}I_{\{x_{k}<0\}}|\s^{1}_{k}(t,x_{k}+x_{k}')-\s^{1}_{k}(t,x_{k}')|^{2}\\
&+ \sum\limits_{k=1}^{m}I_{\{x_{k}<0\}}2x_{k}(b^{1}_{k}(t,\d^{k}x+x')-b^{2}_{k}(t,x')-\displaystyle\int_{E}[\g^{1}_{k}(t,\d^{k}x+x',e)-\g^{2}_{k}(t,x',e)]n(de))\\
&+\sum\limits_{k=1}^{m}I_{\{x_{k}<0\}}\displaystyle\int_{E}[|(x_{k}+\g^{1}_{k}(t,x+x',e)-\g^{1}_{k}(t,\d^{k}x+x',e))^{-}|^{2}-|x_{k}|^{2}\\
&\mbox{ \ \ \ \ \ \ \
 \ \ \ \ \ \ \ \ \ \ \ \ \ \ }-2x_{k}(\g^{1}_{k}(t,x+x',e)-\g^{1}_{k}(t,\d^{k}x+x',e))]n(de)\\&
+\sum\limits_{k=1}^{m}I_{\{x_{k}\geq
0\}}\displaystyle\int_{E}|(\g^{1}_{k}(t,x+x',e)-\g^{1}_{k}(t,\d^{k}x+x',e))^{-}|^{2}n(de)\\
\leq &
\sum\limits_{k=1}^{m}I_{\{x_{k}<0\}}|\s^{1}_{k}(t,x_{k}+x_{k}')-\s^{1}_{k}(t,x_{k}')|^{2}\\
&+\sum\limits_{k=1}^{m}\displaystyle\int_{E}|\g^{1}_{k}(t,x+x',e)-\g^{1}_{k}(t,\d^{k}x+x',e)|^{2}n(de)\\
\leq & (\mu^{2}+\displaystyle\int_{E}\rho^{2}(e)n(de))|x^{-}|^{2}\\
\leq & C^{*}|x^{-}|^{2}.
\end{array}
 \eqno{\Box}
$$

{\bf Remark 3.2.} For the holding of comparison theorem for
multidimensional SDEs with jumps, condition (ii) in Theorem 3.1 is
very natural. $\s^{1}\equiv \s^{2}$ and condition (a) are the
results of that the sign of $dW$ is not always positive or negative;
To consider condition (b) and (c), let us transform our SDEs to the
following forms:
$$\begin{array}{rl}
X^{i}_{s}=&x^{i}+\displaystyle\int_{t}^{s}[b^{i}(r,X^{i}_{r})-\displaystyle\int_{E}\g^{i}(r,X^{i}_{r-},e)n(de)]dr+\displaystyle\int_{t}^{s}\s^{i}(r,X^{i}_{r})dW_{r}\\
&+\displaystyle\int_{t}^{s}\displaystyle\int_{E}\g^{i}(r,X^{i}_{r-},e)N(drde).\end{array}
$$
So condition (b) implies that jumps should occur in the way of
holding the advantage. While condition (c) display the form that the
new drift coefficients should satisfy. In the classical real-valued
case without jumps, we are very familiar with this form.

{\bf Corollary 3.3.} Let $m=1$ and suppose  that
$(b^{i},\s^{i},\g^{i}),i=1,2$ satisfy (A1)and (A2). Then the
following are equivalent:

{\rm(i)} For any $t\in [0,T]$, $x^{1},x^{2}\in R$ such that
$x^{1}\geq x^{2}$, the unique adapted solutions $X^{1}$ and $X^{2}$
in $\mathcal{S}^{2}_{[t,T]}$ to the SDE (3.1) over time interval
$[t,T]$ satisfy:
$$
X^{1}_{s}\geq X^{2}_{s},s\in[t,T],P-a.s.;
$$

{\rm(ii)}For any $t\in [0,T],x\in R$,
$$
\left\{ \begin{array}{l}  \s^{1}(t,x)\equiv\s^{2}(t,x),\\
 b^{1}(t,x)-\displaystyle\int_{E}\g^{1}(t,x,e)n(de)\geq b^{2}(t,x)-\displaystyle\int_{E}\g^{2}(t,x,e)n(de),\\
 x_{1}+\g^{1}(t,x_{1},e)\geq
x_{2}+\g^{2}(t,x_{2},e),\forall x_{1}\geq x_{2}, n(de)-a.s..
\end{array}
\right.  $$

This has already been established in [4].

{\bf Corollary 3.4.} Let $m=1$. When $\g^{1}\equiv\g^{2}\neq 0$ and
suppose that $(b^{i},\s^{i},\g^{1}),i=1,2$ satisfy (A1)and (A2).
Then the following are equivalent:

{\rm(i)} For any $t\in [0,T]$, $x^{1},x^{2}\in R$ such that
$x^{1}\geq x^{2}$, the unique adapted solutions $X^{1}$ and $X^{2}$
in $\mathcal{S}^{2}_{[t,T]}$ to the SDE (3.1) over time interval
$[t,T]$ satisfy:
$$
X^{1}_{s}\geq X^{2}_{s},s\in[t,T],P-a.s.;
$$

{\rm(ii)}For any $t\in [0,T],x\in R$,
$$
\left\{ \begin{array}{l}  \s^{1}(t,x)\equiv\s^{2}(t,x),\\
b^{1}(t,x)\geq b^{2}(t,x),\\
 x_{1}+\g^{1}(t,x_{1},e)\geq
x_{2}+\g^{2}(t,x_{2},e),\forall x_{1}\geq x_{2}, n(de)-a.s..
\end{array}
\right.
$$

{\bf Corollary 3.5.} Let $m=1$. When $\g^{1}\equiv\g^{2}\equiv 0$
and suppose that $(b^{i},\s^{i}),i=1,2$ satisfy (A1) and (A2). Then
the following are equivalent:

{\rm(i)} For any $t\in [0,T]$, $x^{1},x^{2}\in R$ such that
$x^{1}\geq x^{2}$, the unique adapted solutions $X^{1}$ and $X^{2}$
in $\mathcal{S}^{2}_{[t,T]}$ to the SDE (3.1) over time interval
$[t,T]$ satisfy:
$$
X^{1}_{s}\geq X^{2}_{s},s\in[t,T],P-a.s.;
$$

{\rm(ii)}For any $t\in [0,T],x\in R$, $$
\s^{1}(t,x)\equiv\s^{2}(t,x),\mbox{ \ }b^{1}(t,x)\geq b^{2}(t,x). $$

This is the classical result in SDEs without jump.

\vspace {0.2cm}

Although when $\g^{1}\equiv\g^{2}$, it is very convenient for us to
get the comparison theorem. But in fact, it's not necessary for the
holding of comparison theorem. The following is an counter-example,
where
 $\g^{1}$ is not necessarily equal to $\g^{2}$, but the comparison theorem can still hold true.

\vspace {0.2cm}

{\bf Example 3.6.} Let $m=1$. Set
$$b^{i}(t,x)=\int_{E}\g^{i}(t,x,e)n(de),\s^{i}\equiv0,i=1,2.$$
We have the following two SDEs:
$$
X^{i}_{s}=x^{i}+\int_{t}^{s}\int_{E}\g^{i}(r,X^{i}_{r-},e)N(drde).
$$
Then we can immediately see that as long as
$$x_{1}+\g^{1}(t,x_{1},e)\geq
x_{2}+\g^{2}(t,x_{2},e),\forall x_{1}\geq x_{2},n(de)-a.s.,$$ the
comparison theorem: {\rm(i)} in Corollary 3.3 holds true. While in
this case, we only need, for all $t\in [0,T], x\in R$,
$$\g^{1}(t,x,e)\geq
\g^{2}(t,x,e),n(de)-a.s..$$

\subsection{Matrix-valued SDEs}\label{sec:intro}
\qquad Since the first and the second derivatives of the function
$d^{2}_{\mathbb{S}_{+}^{m}}(y), y\in {\mathbb{S}^{m}}$ have been
studied due to Hu and Peng [3], where  $\mathbb{S}^{m}$ is the space
of symmetric real $m\times m$ matrices, and $\mathbb{S}_{+}^{m}$ is
the subspace of $\mathbb{S}^{m}$ containing the nonnegative elements
in $\mathbb{S}^{m}$. So at the end of this paper, we can study the
comparison theorem for matrix-valued SDEs with jumps. Without loss
of generality, we set $d=1$.

Consider the following two SDEs (i=1,2):
$$
X^{i}_{s}=x^{i}+\int_{t}^{s}b^{i}(r,X^{i}_{r})ds+\int_{t}^{s}\s^{i}(r,X^{i}_{r})dW_{r}+\int_{t}^{s}\int_{E}\g^{i}(r,X^{i}_{r-},e)\tilde{N}(drde),
\eqno{(3.3)}
$$
where, for $i=1,2,$
$$
b^i:[0,\infty)\times \mathbb{S}^{m}\rightarrow \mathbb{S}^{m},\mbox{
\ }\s:[0,\infty)\times \mathbb{S}^{m}\rightarrow
\mathbb{S}^{m},\mbox{ \ }\g:[0,\infty)\times \mathbb{S}^{m}\times
R^{l}\rightarrow \mathbb{S}^{m}.
$$

{\bf Theorem 3.7.} Suppose that $b^i,\s^i,\g^i$(i=1,2) satisfy (A1)
and (A2). Then the following are equivalent:

{\rm(i)} For any $t\in [0,T]$, $x^{1},x^{2}\in \mathbb{S}^{m}$ such
that $x^{1}\geq x^{2}$, the unique adapted solutions $X^{1}$ and
$X^{2}$ in $\mathcal{S}^{2}_{[t,T]}(\mathbb{S}^{m})$ to the SDE
(3.3) over time interval $[t,T]$ satisfy:
$$
X^{1}_{s}\geq X^{2}_{s},s\in[t,T],P-a.s.;
$$

{\rm(ii)} $\forall t\in [0,T],\forall (x,x')\in \mathbb{S}^{m}\times
\mathbb{S}^{m}$,
$$\begin{array}{ll}
&-4\langle x^{-},b^1(t,x^++x')-b^2(t,x')\rangle\\
&+\langle
D^{2}d^{2}_{\mathbb{S}_{+}^{m}}(y)(\s^1(t,x+x')-\s^2(t,x')),(\s^1(t,x+x')-\s^2(t,x'))\rangle\\
&+2\displaystyle\int_{E}
[\|(x+\g^1(t,x+x',e)-\g^2(t,x',e))^{-}\|^{2}-\|x^{-}\|^{2}\\
& \ \ \  \ \ \ \ \ \ +2\langle
x^{-},\g^1(t,x+x',e)-\g^2(t,x',e)\rangle]n(de)\\
\leq & C^*\|x^{-}\|^{2},
\end{array}
$$
where $C^{*}\geq 4\mu+\mu^{2}+\int_{E}\rho^{2}(e)n(de)$ is a
constant which does not depend on $t,x,x'$.

{\bf Proof: }We can see from the appendix of [3] that, for any $y\in
\mathbb{S}^{m}$, $y$ has an expression:
$$
y(\l,A)=e^{A}\sum\limits_{i=1}^{m}\l_{i}e_{i}e^{T}_{i}e^{-A},
$$
where $A$ is an antisymmetric real $m\times m$ matrix $(A^{T}=-A)$,
$\l_{i}\in R$, $\{e_{1},e_{2},...,e_{m}\}$ is the standard basis of
$R^{m}$.

If we set
$$
y^{+}(\l,A)=e^{A}\sum\limits_{i=1}^{m}\l^{+}_{i}e_{i}e^{T}_{i}e^{-A},\mbox{
\
}y^{-}(\l,A)=e^{A}\sum\limits_{i=1}^{m}\l^{-}_{i}e_{i}e^{T}_{i}e^{-A}.
$$
Then from [3], we have
$$
d^{2}_{\mathbb{S}_{+}^{m}}(y)=\|y^{-}\|^{2},
\Pi_{\mathbb{S}_{+}^{m}}(y)=y^{+},\mbox{ \ and \ }\nabla
d^{2}_{\mathbb{S}_{+}^{m}}(y)=-2y^{-},
$$
where $\|y\|=(tr(y^{2}))^{\frac{1}{2}}$. This with Lemma 2.3, we can
use the same method as Theorem 3.1 to finish the proof of the
theorem. We omit it.

\end{document}